\documentclass[12pt,reqno]{amsart}
\usepackage{amsfonts,amssymb,hyperref}

\numberwithin{equation}{section}

\begin{document}

\newtheorem{theorem}{Theorem}[section]
\newtheorem{prop}[theorem]{Proposition}
\newtheorem{lemma}[theorem]{Lemma}
\newtheorem{cor}[theorem]{Corollary}

\theoremstyle{definition}
\newtheorem{defn}[theorem]{Definition}
\newtheorem{rmk}[theorem]{Remark}

\newcommand{\C}{{\mathbb C}}
\newcommand\OO{{\mathcal O}}
\newcommand\D{{\mathbb D}}
\newcommand\Z{{\mathbb Z}}
\newcommand\R{{\mathbb R}}

\newcommand\BB{{\mathcal B}}
\newcommand\GG{{ G}}
\newcommand\SSS{{\mathcal S}}

\title[One-relator Sasakian groups]{One-relator Sasakian groups}

\author[I. Biswas]{Indranil Biswas}

\address{School of Mathematics, Tata Institute of Fundamental
Research, Homi Bhabha Road, Mumbai 400005, India}

\email{indranil@math.tifr.res.in}

\author[M. Mj]{Mahan Mj}

\address{School of Mathematics, Tata Institute of Fundamental
Research, Homi Bhabha Road, Mumbai 400005, India}

\email{mahan@math.tifr.res.in}

\subjclass[2000]{57M50, 32Q15, 57M05 (Primary); 14F35, 32J15 (Secondary)}

\date{\today}

\thanks{Both authors were supported by the Department of Atomic Energy, Government of India, 
under project no. 1303/3/2019/R\&D-II/DAE/13820 as also by their respective DST JC Bose 
Fellowships. MM was supported in part by an endowment of the Infosys Foundation, Matrics 
research project grant MTR/2017/000005, CEFIPRA project No. 5801-1.}

\begin{abstract}
We prove that any one-relator group $G$ is the fundamental group of a compact Sasakian manifold 
if and only if $G$ is either finite cyclic or isomorphic to the fundamental group of a compact 
Riemann surface of genus $g \,>\, 0$ with at most one orbifold point of order $n\, \geq\, 1$. We 
also classify all groups of deficiency at least two that are also the fundamental group of some 
compact Sasakian manifold.
\end{abstract}

\maketitle

\section{Introduction}

Fundamental groups of compact Sasakian manifolds are called {\it Sasakian groups}. While the
K\"ahler groups have been studied for long (see \cite{abckt} and references therein), the study 
of their odd-dimensional relatives, viz.\ Sasakian groups, has been started relatively recently 
\cite{chen,bm-higgssas}. In this paper, we extend the results of \cite{bm-onerel,kotschick-def} 
to Sasakian groups by analyzing Sasakian.one-relator groups and Sasakian groups of deficiency 
greater than one.
 
We prove the following generalizing the corresponding statement for K\"ahler groups 
\cite{kotschick-def} (see Theorem \ref{thm-defgenus}):

\begin{theorem}\label{thm-i1}
Let $G$ be a Sasakian group with ${\rm def}(G) \,>\, 1$. Then $G$ must be
an orbifold surface group with genus greater than 1.
\end{theorem}

The following theorem
classifies one-relator Sasakian groups (see Theorem \ref{thm-1rel}):

\begin{theorem}\label{th-onerel-intro}
Let $G$ be an infinite one-relator group. Then $G$ is Sasakian if 
and only if it is isomorphic to 
$$
\langle a_1,\, b_1,\, \cdots ,\, a_i,\, b_i,\, \cdots,\, a_g,\, b_g\, \mid\, 
(\prod_{i=1}^g [a_i,\, b_i])^n 
\rangle\, ,
$$
where $g$ and $n$ are some positive integers.
\end{theorem}

We observe that each of the groups
$$
\langle a_1,\, b_1,\, \cdots ,\, a_g ,\, 
b_g\, \mid\, (\prod_{i=1}^g [a_i,\, b_i])^n \rangle\, ,\,\, \
g,\, n\, >\, 0,
$$
can in fact be realized as the fundamental group of a closed Sasakian manifold as follows. It was
shown in \cite{bm-onerel} that these groups are fundamental groups of smooth projective varieties.
Since the fundamental group of a smooth projective variety is also a Sasakian group
\cite[Proposition 1.2]{chen}, the above groups are also Sasakian groups.

It was shown in \cite{chen} that all finite groups
are Sasakian. Since the only finite one-relator groups are finite cyclic
groups, it follows that
finite one-relator Sasakian groups are precisely the finite cyclic groups.

A useful tool we establish along the way is the following generalization of the
corresponding fact for K\"ahler groups
in \cite[Ch. 4]{abckt} (see Theorem \ref{thm-sasl2}):

\begin{theorem}\label{thm-sasl2intro}
If $G$ is a Sasakian group that satisfies
$\beta_1^2(G)\, >\, 0$, then $G$ is virtually a surface group. 
\end{theorem}

In Theorem \ref{thm-sasl2intro} $\beta_1^2(G)$ denotes the first
$l^2$ betti number of $G$.

\section{Sasakian groups of positive deficiency}

We refer to \cite{BG} for the definition and basic properties of Sasakian manifolds. 
The Sasakian structure on any Sasakian manifold $M$ can be perturbed in
order to make the Sasakian structure quasi-regular \cite{Ru},
\cite[p.~161, Theorem~1.2]{OV}. In
particular, every Sasakian group is the fundamental group of some closed quasi-regular
Sasakian manifold. In view of this, henceforth all Sasakian manifolds considered here
will be assumed to be quasi-regular.

Let $M$ be a quasi-regular closed Sasakian manifold. The group ${\rm U}(1)$ acts on $M$;
let
\begin{equation}\label{al}
\alpha\, :\, M\times \, {\rm U}(1)\, \longrightarrow\, M
\end{equation}
be the action map. The K\"ahler orbifold base $M/{\rm U}(1)$ of $M$ will be denoted by $B$ \cite{Ru},
\cite[p.~208, Theorem~7.1.3]{BG}. Let
\begin{equation}\label{eq-fmb}
f\, :\, M\, \longrightarrow\, M/{\rm U}(1)\,=\, B
\end{equation}
be the quotient map. Let
$$
G\,=\, \pi_1(M)
$$
be the fundamental group. For notational convenience we will omit mentioning the base point
for any fundamental group.

Take a regular point $x_0\, \in\, B$. This means that the action of ${\rm U}(1)$ on
the fiber $f^{-1}(x_0)$ is free.
Let $i\,:\, f^{-1}(x_0) \, \longrightarrow\, M$ be the inclusion map of the regular fiber.
The map ${\rm U}(1)\, \longrightarrow\, f^{-1}(x_0)$ defined by $\lambda\, \longrightarrow\,
x_1\lambda$, where $x_1\, \in\, f^{-1}(x_0)$ is any given point, is a diffeomorphism.
The image $i_\ast \pi_1(f^{-1}(x_0)) \,=\, i_\ast \Z$ is a (finite or infinite) cyclic group
which we shall denote by $C$; to clarify, $C$ may be the trivial group.
Let
$$
Q\,:=\, \pi_1^{orb}(B)
$$
be the orbifold fundamental group of $B$ regarded as the orbifold
quotient of $M$ by the action of ${\rm U}(1)$.
Then there is a short exact sequence:
\begin{equation}\label{eq-sases}
1 \, \longrightarrow\, C \, \longrightarrow\, G\,=\, \pi_1(M) \,
\longrightarrow\, Q \, \longrightarrow\, 1 \,.
\end{equation}

\begin{lemma}[{\cite{BlGo}, \cite{Fu}, \cite{Ta}, \cite{chen}}]\label{lem-b1even}
For the group $G$ in \eqref{eq-sases}, the first betti number $b_1(G)$ is even.
\end{lemma}

Note that from Lemma \ref{lem-b1even} it follows that a nontrivial free group with finitely
many generators can't be a Sasakian group.

The \textit{deficiency} of a finite presented group $\Gamma$ is the maximum of $n-r$ taken over all
possible finite presentations of $\Gamma$, where $n$ and $r$ are the numbers of generators of relations
respectively.

\begin{cor}\label{cor-b1}
If $G$ in \eqref{eq-sases} is of positive deficiency, then $b_1(G) \,\geq\, 2$. If 
$G$ in \eqref{eq-sases} is an infinite one-relator group, then $b_1(G) \,\geq\, 2$.
\end{cor}

\begin{proof}
Since $G$ has positive deficiency, we have $b_1(G) \,\geq\, 1$. From Lemma \ref{lem-b1even} it follows that
$b_1(G) \,\geq\, 2$.

Since infinite one-relator groups have positive deficiency, the second statement follows from the
first statement.
\end{proof}

We shall denote the first $l^2$ betti number of a group $H$ by $\beta_1^2(H)$.
As Kotschick did in \cite{kotschick-def}, we shall make essential use of the following
theorem due to Gaboriau and L\"uck.

\begin{theorem}[{\cite{luck}, \cite{gab}}]\label{thm-luckgab}
If $1 \, \longrightarrow\, N \, \longrightarrow\, H \, \longrightarrow\, Q_1 \, \longrightarrow\, 1$
is a short exact sequence of infinite finitely generated groups with 
$H$ finitely presented, then $\beta_1^2(H) \,= \,0$.
\end{theorem}

\begin{cor}\label{cor-infinitenormal}
Assume that $G$ in \eqref{eq-sases} is of positive deficiency
(for example, it is an infinite one-relator group). 
If $C$ in \eqref{eq-sases} is infinite, then $\beta_1^2(G) \,=\, 0$.
\end{cor}

\begin{proof}
In view of Theorem \ref{thm-luckgab} it suffices to show that $Q$ in \eqref{eq-sases} is infinite.
To prove by contradiction, suppose $Q$ is finite. Then
$G$ has a subgroup $C\,=\,\Z$ of finite index in it. Let $N$ be the finite covering of $M$ such that
$C\,=\, \pi_1(N)\, \subset\, \pi_1(M)\,=\, G$. So $N$ is a quasiregular closed Sasakian manifold
with fundamental group $\Z$. This contradicts Lemma \ref{lem-b1even}.
Hence we conclude that $Q$ is infinite.
\end{proof}

We shall also need the following fundamental inequality:

\begin{lemma}[{\cite{gromov-cras}}]\label{lem-defl2}
Let ${\rm def}(H)$ denote the deficiency of $H$. Then ${\rm def}(H) -1 \,\leq\, \beta_1^2(H)$.
\end{lemma}

Combining Corollary \ref{cor-infinitenormal} and Lemma \ref{lem-defl2} the following is obtained.

\begin{lemma}\label{lem-ruleoutinfinitecyclic}
Assume that $G$ in \eqref{eq-sases} is of deficiency
at least two. Then the cyclic group $C$ in \eqref{eq-sases} is finite.
\end{lemma}

\begin{proof}
Since $G$ is of deficiency at least two,
Lemma \ref{lem-defl2} says that $\beta_1^2(G)\,\geq\, 1$. On the other hand,
if $C$ is infinite, then $\beta_1^2(G)\,=\,0$ by Corollary \ref{cor-infinitenormal}. Hence
the cyclic group $C$ in \eqref{eq-sases} is finite.
\end{proof}

The following lemma is proved in \cite{kotschick-def}.

\begin{lemma}[{\cite[Lemma 3]{kotschick-def}}]\label{lem-lem3kot}
If $C$ in \eqref{eq-sases} is finite, then $Q$ in \eqref{eq-sases}
satisfies the condition $$\beta_1^2(Q) \,=\, p\cdot\beta_1^2(G)\, ,$$
where $p$ is the order of $C$.
\end{lemma}

\begin{lemma}\label{lem-orbfldproj}
Assume that $G$ in \eqref{eq-sases} is of deficiency at least two. Then $Q$
in \eqref{eq-sases} is the fundamental group of a smooth projective orbifold with 
$\beta_1^2(Q) \,>\, 0$.
\end{lemma}

\begin{proof}
Since $M$ in \eqref{eq-fmb} is a closed quasi-regular Sasakian manifold,
the quotient $B$ in \eqref{eq-fmb} is a projective orbifold.

Since $G$ is of deficiency at least two,
Lemma \ref{lem-ruleoutinfinitecyclic} implies that $C$ in \eqref{eq-sases} is finite.
in view of Lemma \ref{lem-defl2}, the given condition, that the 
deficiency of $G$ at least two, also implies that $\beta_1^2(G) \,>\, 0$. 
Hence $\beta_1^2(Q) \,>\, 0$ by Lemma \ref{lem-lem3kot}.
\end{proof}

\section{Sasakian groups with deficiency greater than one}\label{sec-l2}

\subsection{From Sasakian groups to virtual surface groups}

A group $H$ is said to be {\it virtually a surface group} if some 
finite index subgroup of $H$ is the fundamental group of a closed
surface of positive first Betti number.
The following theorem was proved in {\cite{gromov-cras}, \cite{ABR}, \cite[p.~47, Theorem~4.1]{abckt}}
for K\"ahler groups.

\begin{theorem}\label{thm-sasl2}
If $G$ in \eqref{eq-sases} satisfies the condition that
$\beta_1^2(G)\, >\, 0$, then $G$ is virtually a surface group.
\end{theorem}

\begin{proof} Let $C,\, G,\, Q$ be as in \eqref{eq-sases}.
As we saw in the proof of Corollary \ref{cor-infinitenormal}, if $Q$ is finite then
$C$ is finite. On the other hand, if $Q$ is infinite, then Theorem \ref{thm-luckgab}
implies that $C$ is finite. Therefore, we conclude that $C$ is finite.

Hence, by Lemma \ref{lem-lem3kot} we have $\beta_1^2(Q)\,>\,0$.
Note that $Q$ equals the orbifold fundamental group of
the orbifold $B$ in \eqref{eq-fmb}.

Consider the Sasakian metric $g_M$ on the quasi-regular Sasakian manifold $M$ in \eqref{eq-fmb}. 
In the proof of \cite[p.~47, Theorem~4.1]{abckt} substitute $(M,\, g_M)$ in place of the 
K\"ahler metric. Then it is straight-forward to check that all results in Section 4.2 and 
Section 4.3 of \cite{abckt} remain valid.

Let
$$
\psi\, :\, \widetilde{M}\, \longrightarrow\, M
$$
be the universal covering of the quasiregular Sasakian manifold $M$ in \eqref{eq-fmb}. Consider the action of
the fundamental group $\pi_1(M)$ on $\widetilde{M}$. Set
$$
\widehat{M}\, :=\, \widetilde{M}/C\, ,
$$
where $C\, \subset\, \pi_1(M)\,=\, G$ is the subgroup in \eqref{eq-sases}. Let
\begin{equation}\label{p0}
p_0\, :\, \widehat{M}\, \longrightarrow\, M
\end{equation}
be the natural projection.

It can be shown that the action of ${\rm U}(1)$ on $M$ in \eqref{al}
canonically lifts to an action of ${\rm U}(1)$ on $\widehat{M}$. To prove this, let
$$\chi\,=\, \alpha\circ (p_0\times {\rm Id}_{{\rm U}(1)})\, :\, \widehat{M}\times {\rm U}(1)\,
\longrightarrow\, M$$
be the composition, where $\alpha$ is the map in \eqref{al}. The image of the homomorphism
$$
\chi_*\, :\, \pi_1(\widehat{M}\times {\rm U}(1))\, \longrightarrow\, \pi_1(M)
$$
is clearly the subgroup $C\, \subset\, G\,=\, \pi_1(M)$. From this it follows that the action of
${\rm U}(1)$ on $M$ canonically lifts to an action of ${\rm U}(1)$ on $\widehat{M}$.
Let
\begin{equation}\label{vp}
\varphi\, :\, \widehat{M}\, \longrightarrow\,\widehat{B}\, :=\, \widehat{M}/{\rm U}(1)
\end{equation}
be the orbifold quotient for this action of ${\rm U}(1)$ on $\widehat{M}$.

So we have a commutative diagram
\begin{equation}\label{ed}
\begin{matrix}
\widehat{M} & \stackrel{p_0}{\longrightarrow} & M\\
\,\,\,\,\Big\downarrow\varphi && \,\,\Big\downarrow f\\
\widehat{B} & \stackrel{q}{\longrightarrow} & B
\end{matrix}
\end{equation}
where $p_0$, $\varphi$ and $f$ are the maps in \eqref{p0}, \eqref{vp} and
\eqref{eq-fmb} respectively. The map $q$ in \eqref{ed} is an \'etale Galois covering, in the orbifold
category, with Galois group $Q\,=\, G/C$.

Following the proof of Theorem 4.14 in \cite[p.~53, Section~4.4]{abckt} we see that
\begin{itemize}
\item there is a proper holomorphic map to the unit disk
\begin{equation}\label{h}
h\, :\, \widehat{B}\, \longrightarrow\, \mathbb{D}^2\,:=\, \{z\, \in\, {\mathbb C}\, \mid\,\
|z|\, <\, 1\}\, ,
\end{equation}

\item and a homomorphism $\rho\, :\, Q\, \longrightarrow\, \text{Aut}(\mathbb{D}^2)\,=\,
\text{PSL}(2,{\mathbb R})$,
\end{itemize}
such that
\begin{enumerate}
\item the fibers of $h$ are connected,

\item $h$ is $Q$-equivariant for the action of $Q$ on $\mathbb{D}^2$ given by $\rho$ and the
Galois action of $Q$ on $\widehat{B}$, and

\item the homomorphism $h^*\, :\, {\mathcal H}^1_{(2)}(\mathbb{D}^2)\, \longrightarrow\,
{\mathcal H}^1_{(2)}(\widehat{B})$ corresponding to $h$ is an isomorphism.
\end{enumerate}

Consider the complex one-dimensional orbifold
$\OO\,=\, {\mathbb D}^2/\rho(Q)$ with orbifold fundamental group $Q_1\,=\,\rho(Q)$.
Since $h$ in \eqref{h} is $Q$--equivariant, we conclude that
\begin{enumerate}
\item $h$ induces a holomorphic (orbifold) map $h_1\,:\, B\, \longrightarrow\, \OO$
inducing a surjective homomorphism $h_{1\ast}\,:\, Q\, \longrightarrow\, Q_1$ of
orbifold fundamental groups, and

\item the fibers of $h_1$ are compact, connected, because $h$ is proper with connected fibers.
\end{enumerate}

Consequently, there exists a short exact sequence
\begin{equation}\label{K}
1 \, \longrightarrow\, K\, \longrightarrow\, Q
\, \longrightarrow\, Q_1\, \longrightarrow\, 1
\end{equation}
given by $\rho\,=\,h_{1\ast}$. It was observed at the beginning of the proof that
$\beta_1^2(Q) \,>\, 0$. Also, note that $Q_1$ is an infinite group, because $B$ being compact
does note have any nonconstant map to a finite quotient of $\mathbb{D}^2$.
In view of these, by Theorem \ref{thm-luckgab}, the group $K$ in \eqref{K} must be finite.

Since $\OO$ is a compact complex one-dimensional orbifold (it is compact because
$B$ has a nonconstant maps to it), we conclude that $Q_1$ is virtually a surface group. Passing
to a further finite index subgroup $Q_2$
of $Q_1$ if necessary, we may assume that $Q_2$ is a surface group 
acting trivially on the finite normal subgroup $K$ in \eqref{K}. Hence $Q$
is virtually a surface group.

Finally, since $C$ is finite, as pointed at the beginning of the proof, the same argument applied to
$G$ now shows that $G$ is virtually a surface group.
\end{proof}

By Theorem \ref{thm-sasl2} and Lemma \ref{lem-defl2}, we immediately have:

\begin{prop}\label{prop-deflargesas}
If $G$ in \eqref{eq-sases} satisfies the condition that ${\rm def}(G)\, >\, 1$, then $G$ is 
virtually a surface group.
\end{prop}

\subsection{From virtual surface groups to orbifold groups}

In this subsection, we classify virtual surface groups of positive deficiency. Let $\mathcal{G}$ be
a virtual surface group, equivalently, there exists an exact sequence:
\begin{equation}\label{eq-vsurface}
1 \, \longrightarrow\, \pi_1(S) \, \longrightarrow\, \mathcal{G}\, \longrightarrow\, F \, \longrightarrow\, 1\, ,
\end{equation}
where $S$ is a compact surface of positive genus, and $F$ is a finite group
acting by automorphisms on $\pi_1(S)$.

\begin{theorem}\label{thm-defgenus}
Assume that $G$ in \eqref{eq-sases} satisfies the condition ${\rm def}(G) \,>\, 1$.
Then $G$ must be an orbifold surface group with genus greater than 1.
\end{theorem}

\begin{proof} By Proposition \ref{prop-deflargesas}, $G$ is a virtual surface group.

\noindent {\bf Euclidean case:}\,
Suppose that the genus of the surface $S$ in \eqref{eq-vsurface} for
$\mathcal{G}\,=\, G$ is one. Then there exists a 2-dimensional
crystallographic group $H$ such that $$G \,=\, H \times F_1\, ,$$ where $F_1$ is the subgroup of
$F$ in \eqref{eq-vsurface} acting
trivially on $\pi_1(S)$. The deficiency of any crystallographic group other than $\Z\times\Z$, and the Klein
bottle group, is non-negative. The deficiency of any finite group is also non-negative. Since 
$\Z\times\Z$ and the Klein bottle group 
have deficiency exactly one, and need at least two generators, it follows that if $F_1$ is non-trivial,
then ${\rm def}(G)\, \geq\, 0$. Thus, $F_1$ is trivial, and we have one of the two possibilities
for $G$:
\begin{enumerate}
\item $\Z \times \Z$, and

\item the Klein bottle group.
\end{enumerate}

Since the Klein bottle group has first betti number one, it cannot be Sasakian (see Lemma \ref{lem-b1even}). 
Finally, $\Z\times\Z$ has deficiency equal to one and is ruled out by the hypothesis.

\noindent {\bf Hyperbolic case:}\, Next, suppose that the genus of $S$ in \eqref{eq-vsurface} is 
greater than one. By the Nielsen realization theorem \cite{kerckhoff-nr}, there exists a 
hyperbolic orbifold $\OO$ with $\pi_1(\OO) \,=\, H$ such that $$G \,=\, H \times B_1$$ with 
$B_1$ finite. As in the genus one case, it follows that $B_1$ is trivial. If $\OO$ is 
non-orientable, then $b_1(\OO)$ is odd, forcing $\OO$ to be an orientable hyperbolic orbifold 
(using Lemma \ref{lem-b1even} again). Further, from the given condition that the deficiency is 
greater than one it follows that the genus of $\OO$ is greater than one.
\end{proof}

\section{One relator groups}

\subsection{One relator Sasakian groups are projective orbifold groups}

Murasugi has described in detail the centers of one-relator groups.

\begin{theorem}[{\cite[Theorems 1, 2]{mu}}]\label{center} Let 
$\mathcal{G}\,=\, \langle x_1,\, \cdots, \, x_k\, \mid\, w^n\rangle$ be a 
one-relator group with $n\,\geq\,1$. If $k\,\geq\, 3$, then the center 
$Z(\mathcal{G})$ of $\mathcal{G}$ is trivial. If $\mathcal{G}$ is not abelian with $k\,=\,2$,
and $Z(\mathcal{G})$ is not trivial, then $Z(\mathcal{G})$ is infinite cyclic.
\end{theorem} 

\begin{theorem}[{\cite[p.~219]{ks}}]\label{ks} Let $\mathcal{G}$ be a 
one-relator group having a (nontrivial) finitely presented normal 
subgroup $H$ of infinite index. Then $\mathcal{G}$ is torsion-free and has 
two generators. Further, $\mathcal{G}$ is an infinite cyclic or infinite dihedral 
extension of a finitely generated free group $N$ (meaning $N\,\subset\,\mathcal{G}$) satisfying the following: 
\begin{itemize}
\item $H \,\subset\, N$ if $H$ is not cyclic, and
		
\item $H\bigcap N$ is trivial if $H$ is cyclic.
\end{itemize}
\end{theorem}

\begin{prop}\label{prop-ctrivial}
Assume that $G$ in \eqref{eq-sases} is an infinite one relator group.
Then
\begin{itemize}
\item $C$ in \eqref{eq-sases} is trivial, and

\item $G\,=\, Q$ is a projective orbifold group.
\end{itemize}
\end{prop}

\begin{proof}
By Theorem \ref{center}, $C$ is trivial if ${\rm def}(G)\,>\,1$. In that case, $G\,=\,Q$ is a 
projective orbifold group.

On the other hand, if $G$ is abelian, then $G\, =\, \Z\times\Z$ and $G$ is a projective orbifold 
group, as it is the fundamental group of an elliptic curve.

Therefore, assume that ${\rm def}(G)\,=\,1$ and $G$ is non-abelian.

Since ${\rm def}(G)\,=\,1$, the minimum number of generators of the one relator group $G$, which 
we will denote by $k$, is two.

Since $C$ is contained in the center $Z(G)$ of $G$, if $Z(G)$ is trivial, then the proposition follows.
So we assume that $Z(G)$ is non-trivial.

By Theorem \ref{center}, $C\,=\,\Z$. Hence, by Theorem
\ref{ks}, $G\,=\, N \times \Z$, where $N\,=\,F_r$ is free of rank $r$ greater than one; we note that
if $N$ is free of rank one, then $G\,=\, \Z \times \Z$ and $C \,\neq\, \Z$.

Next, since $G$ is one relator, we have $b_2(G) \,\leq\, 1$. On the other hand, $$b_2(N \times 
\Z) \,=\, r \,>\, 1\, ,$$ so we get a contradiction. Thus, if $k\,=\,2$, then $Z(G)$ is trivial, 
forcing $C$ to be trivial. Consequently, $G\,=\,Q$ is a projective orbifold group.
\end{proof}

We shall need the following result.

\begin{prop}[{\cite[Theorem 1]{fks}}]\label{fks} Let
$\mathcal{G}\,=\, \langle x_1,\, \cdots, \, x_k\, \mid\, w^n\rangle$ be a one-relator group with 
$n\, \geq \, 1$. If $k\,=\,1$, then $\GG$ is torsion-free. Else, every torsion element in 
$\mathcal{G}$ is conjugate to a power of $w$, and the subgroup generated by torsion elements in 
$\mathcal{G}$ is the free product of the conjugates of $w$.
\end{prop}

\subsection{From orbifold projective groups to projective groups}\label{sec-orbfldproj2proj}

In this subsection we assume that
$G$ in \eqref{eq-sases} is a one-relator Sasakian group
$$G \,=\, \langle x_1,\, \cdots, \, x_k\,\mid\, w^n \rangle\, .$$ 
If ${\rm def}(G)\, >\, 1$, Theorem \ref{thm-defgenus} shows that
$G\,=\,\pi_1(\OO)$, where $\OO$ is a hyperbolic orbifold of genus greater than one.

Therefore, we assume that ${\rm def}(G)\, =\, 1$, so $k\,=\,2$.

Since $k\,=\,2$, 
Proposition \ref{prop-ctrivial} further shows that
\begin{equation}\label{gb}
G\,=\,\pi_1(B)\, ,
\end{equation}
where $B$ is the projective orbifold in \eqref{eq-fmb}.

It can be shown that $\dim_\C B\,>\,1$. Indeed, if $\dim_\C B\,=\,1$, then $M$ in
\eqref{eq-fmb} is a 3-dimensional Seifert-fibered
manifold. Since $G\,=\,\pi_1(M)$ is infinite, it follows that the center of $G$ is infinite cyclic
\cite[Ch.~12]{hempel-book}. But this contradicts Proposition \ref{prop-ctrivial}.
So we conclude that $\dim_\C B\,>\,1$.

\noindent {\bf Structure of singularities:}\, We analyze the locus of the points of the orbifold 
$B$ in \eqref{eq-fmb} with nontrivial inertia group. Such examples exist even when, at the 
global level, the Sasakian manifold is simply connected \cite{CMST,MST,MT}. Let $\mathcal S$ be 
a connected component of this locus of $B$. Let $U_{\mathcal S}$ denote a regular neighborhood 
of $\mathcal S$ in $B$.

\begin{lemma}\label{lem-pi1ofsingnbhd} Let $2m$ denote the real dimension of $B$.
There exists a lens space
$L$ of dimension $2l-1$, with $l \,\geq \,2$, such that the topological space for the
orbifold $U_{\mathcal S}$ is homeomorphic to $cL \times
{\D}_{2m-2l}$, where $cL$ denotes the cone on $L$ and ${\D}_{2m-2l}$ denotes a ball of real dimension $2m-2l$.
\end{lemma}

\begin{proof}
We recall $B\,=\, M/{\rm U}(1)$ with $M$ being a smooth manifold. All the isotropy subgroups for the
action of ${\rm U}(1)$ on $M$ are finite cyclic groups. The local structure of $B$ follows from this
(see \cite[p.~146, Theorem 4.7.7]{BG}).
\end{proof}

Let $i_{\mathcal S}\,:\, U_{\mathcal S} \,\longrightarrow\, B$ be the inclusion map
between orbifolds, and let
$i_{{\mathcal S},\ast}$ denote the homomorphism between orbifold fundamental groups induced
by $i_{\mathcal S}$.

\begin{lemma}\label{lem-alltrivial}
Suppose $i_{{\mathcal S},\ast}(\pi_1^{orb} (U_{\mathcal S}))$ is trivial for all components
$\mathcal S$. Then the group $G$ is projective. In particular, if $G$ is torsion-free, it is projective.
\end{lemma}

\begin{proof}
We may blow-up the orbifold $B$ and obtain another orbifold
$$
\varpi\, :\, B'\, \longrightarrow\, B
$$
such that the underlying topological space for $B'$ is a smooth projective variety. We note that 
the underlying topological space for an orbifold $B''$ is a smooth projective variety if the 
locus in $B''$ of points with nontrivial inertia is a normal crossing divisor. In view of the 
given condition, we may choose $B'$ such that $\pi_1(B)$ coincides with the fundamental group of 
the underlying topological space for $B'$. But the underlying topological space for $B'$ is a 
smooth projective variety, so $\pi_1(B)$ is a projective group.

If $G$ is torsion-free, then $i_{{\mathcal S},\ast}(\pi_1^{orb} (U_{\mathcal S}))$ must be 
trivial for all components $\mathcal S$. The conclusion of the lemma follows.
\end{proof}

\begin{lemma}\label{lem-atmostone}
Suppose $i_{{\mathcal S},\ast}(\pi_1^{orb} (U_{\mathcal S}))$ is non-trivial for some ${\mathcal 
S}$. Then any torsion element of $G$ has a non-trivial power that is conjugate to an element of 
$i_{{\mathcal S},\ast}(\pi_1^{orb} (U_{\mathcal S}))$.
\end{lemma}

\begin{proof}
This is an immediate consequence of Proposition \ref{fks}.
\end{proof}

We next analyze the case where the underlying topological space for $B$ is smooth.

\begin{lemma}\label{lem-smoothbase}
Suppose that the underlying topological space for the projective orbifold $B$ is smooth. Then 
$G$ is isomorphic to the fundamental group of a (real) two-dimensional compact orbifold $V_1$ 
with at most one cone-point. Furthermore, $V_1$ is orientable.
\end{lemma}

\begin{proof}
Let $$
G\,=\, \langle x_1\, ,\cdots \, ,x_k\, \mid\, w^n\rangle\, ,
~\, n\, >\, 1\, .
$$
Then, replacing the orbifold $B$ by its underlying topological space $\BB$, we note by
Proposition \ref{fks} and Lemma
\ref{lem-atmostone}, that $\pi_1(\BB)$ is the quotient of $\pi_1(B)$ by some powers
of $w$. Hence $G_1\,=\,\pi_1(\BB)$ is of the form
$$
G_1\,=\, \langle x_1\, ,\cdots \, ,x_k\, \mid\, w^r\rangle\, ,
~\, r\, \geq\, 1\, .
$$

Since $\BB$ is a smooth projective variety, $G_1$ is a one-relator K\"ahler group. Hence, by the 
classification of one-relator K\"ahler groups \cite{bm-onerel} (see also \cite{kotschick-def}), 
$G_1$ is isomorphic to the fundamental group of a (real) two-dimensional compact orbifold $V$ 
with at most one cone-point $y_0$. If $r\,=\,1$, then $G_1$ is torsion-free. Else the loop that goes 
around $y_0$ represents the conjugacy class of $w$. Replacing $r$ by $n$, we note that $G$ is 
isomorphic to the fundamental group of a (real) two-dimensional compact orbifold $V_1$, where 
$V_1$ differs from $V$ only in the order of the cone-point $y_0$. Thus, $G$ is isomorphic to the 
fundamental group of a (real) two-dimensional compact orbifold $V$ with at most one cone-point. 
\end{proof}

Lemma \ref{lem-pi1ofsingnbhd} furnishes the local structure of singularities of $B$. We now 
describe the blowup. We first note that a lens space $cL$ on a lens space $L$ of dimension $2k-1 
$ with fundamental group $\Z/q\Z$ can be realized as the topological boundary of the total space 
$E_q$ of a twisted line bundle (more precisely, the $q-$th tensor power of the tautological line 
bundle) over $\C P^{k-1}$. Hence, the blow up of $cL$ may be obtained by removing an open 
neighborhood of the cone point $y_0 \,\in\, cL$ and attaching a copy of $E_q$ along the resulting 
boundary $L$. Let $BU(L)$ denote the blowup of $cL$, and let $ED(L) \,\subset\, BU(L)$ be the 
exceptional divisor, which is homeomorphic to $\C P^{k-1}$. Then we note the following two
properties of $BU(L)$:

\begin{rmk}\label{rmk-blowuplens}\mbox{}
\begin{enumerate}
\item If $cL$ has an orbifold cone-point of order $q$ then $BU(L)$ is also an orbifold so that the exceptional 
divisor $ED(L)$ is an orbifold locus of ramification order $q$.

\item The underlying topological space of $BU(L)$ is a smooth manifold.
\end{enumerate}
\end{rmk}

Next, we describe the local structure of a singularity of $B$ as obtained in Lemma 
\ref{lem-pi1ofsingnbhd} along with the ${\rm U}(1)$ bundle over it. Let $E_\SSS$ denote the circle 
bundle over $U_\SSS$ induced by the inclusion of $U_\SSS$ into $B$. Thus, we have
$E_\SSS\, \subset\, M$. 
Then we have the local structure of the Sasakian manifold $M$ given by the following commutative 
diagram:
\begin{equation}\label{locstr}
\begin{matrix}
D^{2k} \times U(1) & {\longrightarrow} & E_\SSS\\
\,\,\,\,\Big\downarrow && \,\,\Big\downarrow \\
\, \quad \, D^{2k} & {\longrightarrow} & U_\SSS = cL
\end{matrix}
\end{equation}

\begin{theorem}\label{thm-1rel} Let
$$
G\,=\, \langle x_1\, ,\cdots \, ,x_n\, \mid\, w^k\rangle\, , ~\ k\, >\, 1\, , 
$$
be a one-relator Sasakian group. Then $G$ is isomorphic to the fundamental group of a (real) 
two-dimensional compact orbifold $V$ with at most one cone-point. Further, the underlying 
manifold of $V$ is orientable.
\end{theorem}

\begin{proof}
As discussed at the beginning of Section \ref{sec-orbfldproj2proj}, it remains only to deal 
with the case that $n\,=\,2$. If $k\,=\,1$, then $G$ is torsion-free by Theorem \ref{fks}, and 
the hypothesis of Lemma \ref{lem-alltrivial} is satisfied, forcing $G$ to be the fundamental
group of a smooth projective variety. The main theorem of \cite{bm-onerel} now furnishes
the result.

Else, the local structure of singularities is given by Lemma \ref{lem-pi1ofsingnbhd}.
Blowing up and using Remark \ref{rmk-blowuplens}, we obtain an orbifold $B_1$ satisfying the
following:
\begin{enumerate}
\item The underlying topological space $\BB_1$ of $B_1$ is a smooth projective variety.

\item The orbifold locus of $B_1$ is a divisor.

\item Any loop that goes around a component of the divisor represents an element of $\pi_1^{orb} 
(B) \,=\,G$ that is conjugate to $w^r$ for some $r$ that divides $k$.
\end{enumerate} 

Hence, by Lemma \ref{lem-smoothbase}, $G_1\,=\,\pi_1(\BB_1)$ is isomorphic to the fundamental group 
of a (real) two-dimensional compact orbifold $V$ with at most one cone-point. The rest of the 
argument is a replica of the last part of the proof of Lemma \ref{lem-smoothbase} forcing $G$ to 
be isomorphic to the fundamental group of a (real) two-dimensional compact orbifold $V_1$, where 
$V_1$ differs from $V$ only in the order of the cone-point.
\end{proof}

\end{document}